\documentclass[12pt]{article}
\usepackage{amsmath}
\usepackage{amssymb}
\usepackage{amsfonts}
\usepackage{tikz}
\usepackage{authblk}

\usepackage{amsthm}

\newtheorem*{notation*}{Notation}
\newtheorem{theorem}{Theorem}[section]
\newtheorem{definition}[theorem]{Definition}

\newtheorem{lemma}[theorem]{Lemma}

\newtheorem{corollary}[theorem]{Corollary}

\newtheorem*{remark*}{Remark}
\numberwithin{equation}{section}

\DeclareMathOperator{\ex}{ex}

\DeclareMathOperator{\spl}{split}

\def \a {\alpha}

\def \b {\beta}

\def \N {\mathbb{N}}

\def \C {\mathcal{C}}
\def \L {\mathcal{L}}
\def \G {\mathcal{G}_{tree}}
\def \Gf {\mathcal{G}_{forest}}
\def \D {\Delta}

\tikzstyle{cir} = [draw, circle, minimum height= 20 mm]

\title{Exact results on generalized Erd\H{o}s-Gallai problems}
\author[1]{Debsoumya Chakraborti\thanks{Supported by the Institute for Basic Science (IBS-R029-C1), and the European Research Council (ERC) under the European Union Horizon 2020 research and innovation programme (grant agreement No. 947978).}}
\author[2]{Da Qi Chen\thanks{This material is based upon work supported by the Air Force Office of Scientific Research under award number FA9550-20-1-0080.}}
\affil[1]{\small Mathematics Institute, University of Warwick, Coventry, UK}
\affil[2]{\small University of Virginia, Charlottesville, VA, USA}
\affil[ ]{\small Email:
\texttt{debsoumya.chakraborti@warwick.ac.uk},
\texttt{wny7gj@virginia.edu}}

\begin{document}
\maketitle
\begin{abstract}
Generalized Tur\'an problems have been a central topic of study in extremal combinatorics throughout the last few decades. One such problem is maximizing the number of cliques of size $t$ in a graph of a fixed order that does not contain any path (or cycle) of length at least a given number. Both of the path-free and cycle-free extremal problems were recently considered and asymptotically solved by Luo. We fully resolve these problems by characterizing all possible extremal graphs. We further extend these results by solving the edge-variant of these problems where the number of edges is fixed instead of the number of vertices. We similarly obtain exact characterization of the extremal graphs for these edge variants. 
\end{abstract}

\section{Introduction}

One of the seminal problems in extremal graph theory is to determine $\ex(n, P_k)$, which denotes the maximum number of edges in an $n$-vertex graph that does not contain a copy of $P_k$ (a path with $k$ vertices). Erd\H{o}s and Gallai showed in \cite{EG1} that $\ex(n,P_k) \le \frac{k-2}{2} \cdot n$. It is easy to see that when $k-1$ divides $n$, the result is tight and the extremal graph is the vertex-disjoint union of $\frac{n}{k-1}$ copies of $K_{k-1}$. Faudree and Schelp \cite{FS75} and independently Kopylov \cite{Ko} later determined $\ex(n,P_k)$ exactly for all values of $n$ and described the extremal graphs.

\begin{theorem} [\cite{FS75,Ko}] \label{erdos} 
Let $n = q(k-1) + r$, $0 \le r \le k-2$, and $k \ge 2$. Then, $$\ex(n, P_k) = q \binom{k-1}{2} + \binom{r}{2}.$$
Moreover, the extremal graphs are  
\begin{itemize}
\item $qK_{k-1}\cup K_r$, vertex disjoint unions of $q$ complete graphs $K_{k-1}$ and a $K_r$, or
\item when $k$ is even and $r$ is either $\frac{k}{2}$ or $\frac{k}{2} - 1$, another extremal graph can be obtained by taking a vertex disjoint union of $t$ copies of $K_{k-1}$ $(0 \le t < q)$ and a copy of $H$, where $H$ is the graph achieved by adding all the edges between a clique with $\frac{k}{2}-1$ vertices and an independent set with $n-t(k-1)-\frac{k}{2}+1$ vertices. 
\end{itemize}
\end{theorem}

Erd\H{o}s and Gallai also showed in \cite{EG1} that any graph with $n$ vertices and more than $\frac{k-1}{2} \cdot (n-1)$ edges must have a cycle of length at least $k$. It is easy to see that when $k-2$ divides $n-1$, the result is tight if one considers the extremal graph obtained by gluing $\frac{n-1}{k-2}$ copies of $K_{k-1}$ at a single vertex. For the convenience of classification of such extremal graphs, we define the following class of graphs.

\begin{definition} \label{def:family}
Let $F_1, F_2, \ldots, F_l$ be a family of connected graphs. Define $\G(F_1, F_2, \ldots, F_l)$ $(\text{analogously} \; \Gf(F_1, F_2, \ldots, F_l))$ to be the family of graphs $G$ satisfying the following two properties. 
\begin{enumerate}
	\item The edges of $G$ can be partitioned into $l$ disjoint sets $E_1, E_2, \dots, E_l$ such that the graph induced by $E_i$, $G(E_i)$, is isomorphic to $F_i$ for $i=1, \ldots, l$. Furthermore, for $i\neq j$, $|V(G(E_i))\cap V(G(E_j))|\le 1$. 
	\item Let $H$ be the $l$-vertex auxiliary graph (with vertex set $\{1, 2, \ldots, l\}$) created from $G$ such that $ij\in E(H)$ if and only if $G(E_i)$ and $G(E_j)$ shares a vertex. Then, $H$ is a tree (a forest). 
\end{enumerate}
We often refer the graphs $F_1, F_2, \ldots, F_l$ as blocks.
\end{definition} 

Following the literature, we denote by $\ex(n,C_{\ge k})$ the maximum number of edges in an $n$-vertex graph that does not contain a cycle of length at least $k$. The exact value of $\ex(n,C_{\ge k})$ was later determined by Woodall \cite{W} and independently by Kopylov \cite{Ko}.

\begin{theorem} [\cite{Ko,W}] \label{kopylov}
Let $n = q(k-2) + r$, $1 \le r \le k-2$, $k \ge 3$, and $q \ge 1$. Then, $$\ex(n, C_{\ge k}) = q \binom{k-1}{2} + \binom{r}{2}.$$
Moreover, the extremal graphs are  
\begin{itemize}
\item $\G(qK_{k-1}, K_r)$, where the blocks are $q$ complete graphs $K_{k-1}$ and a $K_r$, or
\item when $k$ is odd and $r$ is either $\frac{k+1}{2}$ or $\frac{k-1}{2}$, another class of extremal graphs is $\G(tK_{k-1}, H)$ (for any $0 \le t < q$), where $H$ is the graph achieved by adding all the edges between a clique with $\frac{k-1}{2}$ vertices and an independent set with $n-t(k-2)-\frac{k}{2}+1$ vertices. 
\end{itemize}
\end{theorem}

The field of extremal graph theory is very broad and interested readers can peruse the excellent survey \cite{FS} by F\"uredi and Simonovits for more related problems. In the recent decade, there has been a growing interest in studying the generalized variant where one maximizes the number of copies of a fixed graph instead of the number of edges (see, e.g., \cite{AS,BG,FO,Gr,GL,HH}). Following the notation in \cite{AS}, for two graphs $T$ and $H$, the generalized extremal function $\ex(n,T,H)$ is defined to be the maximum number of copies of $T$ in an $n$-vertex $H$-free graph. Note that for $T = K_2$, we have the standard extremal function where we maximize the number of edges in an $n$-vertex $H$-free graph, i.e., $\ex(n,K_2,H) = \ex(n,H)$. 

For the generalized extremal problems, there exist substantial interests in the case when the forbidden graph $H$ is a tree. In particular, when $H$ is a star (e.g. $K_{1,k}$), after significant progress made in \cite{ACM,AM,CR,CR1,EG,G,GLS,LM}, Chase \cite{C} recently computed $ex(n,K_t,K_{1,k})$ with the classification of all the extremal graphs. As for the case when the forbidden graph is a path, Luo asymptotically answered the generalized versions of Theorems~\ref{erdos} and~\ref{kopylov} in \cite{Luo}, where they proved the following couple of statements as corollaries of more general results.

\begin{theorem} [\cite{Luo}] \label{Luo1}
The extremal function $\ex(n, K_t, P_k) \le \frac{n}{k-1} \binom{k-1}{t}$. 
\end{theorem}

\begin{theorem} [\cite{Luo}] \label{Luo2}
The extremal function $\ex(n, K_t, C_{\ge k}) \le \frac{n-1}{k-2} \binom{k-1}{t}$. 
\end{theorem}

Luo's results were useful in investigating certain Tur\'an-type problems in hypergraphs (see, e.g., \cite{GMSTV}). Notably, Ning and Peng provided some extensions and applications of Theorems~\ref{Luo1} and~\ref{Luo2} in \cite{NP}. In this paper, we strengthen Theorems~\ref{Luo1} and~\ref{Luo2} by proving exact results on $\ex(n, K_t, P_k)$ and $\ex(n, K_t, C_{\ge k})$ with the classification of all extremal graphs. For the convenience of writing, we denote the number of cliques of order $t$ in a graph $G$ by~$k_t(G)$.

\begin{theorem} \label{path}
For any positive integers $n$ and $3 \le t < k$, if $G$ is a $P_k$-free graph on $n$ vertices, then $k_t(G) \le k_t(qK_{k-1} \cup K_r)$, where $n = q(k-1) + r$, $0 \le r \le k-2$. Moreover, $qK_{k-1} \cup K_r$ is the unique graph that achieves equality when $t \le r$. If $t > r$, then $G$ is an extremal graph if and only if $G$ is isomorphic to $qK_{k-1} \cup L$, where $L$ is any graph on $r$ vertices.
\end{theorem}

Note that if $t \ge k$ in Theorem~\ref{path}, then $\ex(n,K_t,P_k) = 0$ because any graph $G$ with a copy of $K_t$ contains a copy of $P_k$. 

\begin{theorem} \label{cycle}
For any positive integers $n$ and $3 \le t < k$, if $G$ is a $C_{\ge k}$-free graph on $n$ vertices, then $k_t(G) \le k_t(qK_{k-1} \cup K_r)$, where $n = q(k-2) + r$, $1 \le r \le k-2$. Furthermore, $G$ is an extremal graph if and only if $G$ is in:
	\begin{itemize}
	\item $G\in \G(qK_{k-1}, K_r)$ when $t\le r$, or
	\item $\mathcal{G}'$ when $t>r$ where $\mathcal{G}'=\{G: |V(G)| = n, G\in \mathcal{G}_{forest}(qK_{k-1}, L) \text{ for some graph } L\}$.
	\end{itemize}
\end{theorem}

\begin{remark*}
Note that we have defined $\Gf$ only for a family of connected graphs. However, observe that $L$ can be a disconnected graph in Theorem~\ref{cycle} when $r < t$. Thus, we extend the definition of $\Gf$ as follows:  let $F_1, F_2, \ldots, F_l$ be an arbitrary family of graphs, and let $F_{i,1}, \ldots, F_{i,j_i}$ be the connected components of $F_i$ for all $1 \le i \le l$. Then, define $\Gf(F_1, F_2, \ldots, F_ l)$ to be $\Gf(F_{1,1}, F_{1,2}, \ldots, F_{1,j_1}, F_{2,1}, \ldots, F_{2,j_2}, \ldots, F_{l,j_l})$. 
\end{remark*}

Another popular variant of these classical extremal problems is to study the edge analogues where the number of edges is fixed instead of the number of vertices (e.g. \cite{CC,F,FFK,KR,KR1}). Motivated by this line of work, we study the edge analogues of the problems in Theorems~\ref{path} and~\ref{cycle} and classify all their extremal graphs. Note that adding or deleting isolated vertices from a graph does not change the number of edges nor the number of $K_t$'s for $t \ge 3$. Thus, in this paper, when we discuss the edge variant, two graphs are considered to be equivalent if they are isomorphic after deleting all the isolated vertices. 

In order to describe the extremal graphs, we need a couple of definitions that are similar to those in \cite{CC}. First, we introduce the notions of colex (colexicographic) order and colex graphs. Colex order on the finite subsets of the natural number set $\mathbb{N}$ is defined as the following: for $A, B \subseteq \mathbb{N}$, we have that $A \prec B$ if and only if $\max((A \setminus B) \cup (B \setminus A)) \in B$. The colex graph $L_m$ on $m$ edges is defined as the graph with the vertex set $\mathbb{N}$ and edges are the first $m$ sets of size $2$ in colex order. Note when $m=\binom{r}{2}+s$ where $0\le s< r$, then $L_m$ is the graph containing a clique of order $r$ and an additional vertex adjacent to $s$ vertices of the clique. The celebrated Kruskal-Katona theorem implies the following.

\begin{theorem} [\cite{Ka,Kr}] \label{maxkt}
For any positive integers $t$, $m$, and any graph $G$ on $m$ edges, we have that $k_t(G) \le k_t(L_m)$. Moreover, $L_m$ is the unique graph satisfying the equality when $s \ge t-1$, where $m = \binom{r}{2} + s$ with $0 \le s < r$.
\end{theorem}

\begin{remark*}
For Theorem~\ref{maxkt}, if $r \ge t$ and $s < t-1$, a graph $G$ satisfies the equality if and only if $G$ is an $m$-edge graph which contains $K_r$ as a subgraph. If $r < t$, then any graph with $m$ edges satisfies the equality because no graph on $m$ edges has a copy of $K_t$.
\end{remark*} 

Next, for notational convenience, we define the following class of graphs.

\begin{definition}
For $m=0$, let $\L_{t,k}(m)$ $(\text{analogously} \; \C_{t,k}(m))$ be the family of empty graph, and for $0 < m \le \binom{k-1}{2}$, call $\L_{t,k}(m)$ $(\C_{t,k}(m))$ to be the following family of graphs, where $m = \binom{r}{2} + s$ with $0 \le s < r$.
\begin{itemize}
\item If $s \ge t-1$, then $\L_{t,k}(m)$ ($\mathcal{C}_{t, k}(m)$) contains just the colex graph $L_m$.
\item If $r \ge t$ and $s < t-1$, then $\L_{t,k}(m)$ ($\mathcal{C}_{t, k}(m)$) contains not only $L_m$, but also all $m$-edge $P_k$-free ($C_{\ge k}$-free) graphs with $K_r$ as a subgraph.
\item If $r < t$, then $\L_{t,k}(m)$ ($\mathcal{C}_{t, k}(m)$) contains not only $L_m$, but also all $m$-edge $P_k$-free ($C_{\ge k}$-free) graphs.
\end{itemize} 
\end{definition}

Next, we state our result on the edge analogue of Theorem~\ref{cycle}.

\begin{theorem} \label{edge}
For any $3 \le t < k$, if $G$ is a $C_{\ge k}$-free graph with $m$ edges, then $k_t(G) \le k_t(qK_{k-1} \cup L_b)$, where $m = q\binom{k-1}{2} + b$ and $0 \le b < \binom{k-1}{2}$. Moreover, $G$ is an extremal graph if and only if $G$ is isomorphic to some graph in $\Gf(qK_{k-1}, L)$ for some $L \in \C_{t,k}(b)$. 
\end{theorem}

As a simple corollary of this theorem, we obtain the edge analogue of Theorem~\ref{path}.

\begin{corollary} \label{edgecor}
For any $3 \le t < k$, if $G$ is a $P_k$-free graph with $m$ edges, then $k_t(G) \le k_t(qK_{k-1} \cup L_b)$, where $m = q\binom{k-1}{2} + b$ and $0 \le b < \binom{k-1}{2}$. Moreover, $G$ is an extremal graph if and only if $G$ is isomorphic to $qK_{k-1} \cup L$ for some $L \in \L_{t,k}(b)$. 
\end{corollary}

The rest of this paper is organized as follows. In Section~\ref{sec:vertex}, we prove Theorem~\ref{path}. In Section~\ref{sec:newedge}, we prove the edge-variant, Theorem~\ref{edge}, by developing and using relevant structural results. We, then, use Theorem~\ref{edge} to prove Theorem~\ref{cycle} and Corollary~\ref{edgecor} in Section~\ref{sec:cor}. We end with a few concluding remarks in Section~5. 

\section{Proof of Theorem~\ref{path}}
\label{sec:vertex}

In this section, our goal is to prove Theorem~\ref{path}. The strategy is to first assume that $G$ is a minimum counter-example. In particular, we assume that there exists a minimum $t\ge 3$ and a $P_k$-free graph $G$ with minimum order $n=q(k-1)+r$  such that $k_t(G) \ge k_t(qK_{k-1}\cup K_r)$ and $G$ is not one of the extremal graphs described in Theorem~\ref{path}. We further assume that the number of copies of $K_t$ in any other graph with the same number of vertices is at most $k_t(G)$. We start by proving some structural results about this minimum counter-example $G$ and lastly show that such an example does not exist. 

Note that if $n \le k-1$, it is easy to check that $G = K_n$ maximizes $k_t(G)$. Moreover, if $t\le n \le k-1$, then $K_n$ is the unique extremal graph. Also, if $n < t$, then any graph on $n$ vertices serves as an extremal graph because no graph on $n$ vertices contains a copy of $K_t$. Thus, for the rest of this section, we may assume that $n \ge k$. We use the following simple fact proved in \cite{CC}.

\begin{lemma} [Lemma 2.3 in \cite{CC}] \label{easy:convex}
Let $t$, $w$, $x$, $y$, and $z$ be non-negative integers such that $t \ge 2$, $x+w=y+z$, $x\ge y$, $x\ge z$, and $x\ge t$. Then, $\binom{x}{t}+\binom{w}{t} \ge \binom{y}{t}+\binom{z}{t}$. Moreover, the inequality is strict if $x> y$ and $x>z$. 
\end{lemma}

\begin{lemma} \label{con}
If $G$ is a minimum counter-example, then $G$ is connected.
\end{lemma}

\begin{proof}
Suppose for the sake of contradiction that $G$ is not connected. First, we show that $G$ does not contain a clique of order $k-1$. Suppose $G$ does contain one such clique; note that due to the $P_k$-free condition, the clique is disjoint from the rest of the graph. Then, by removing this clique, we obtain a graph $G'$ with $(q-1)(k-1)+r$ vertices with at least $(q-1)\binom{k-1}{t}+\binom{r}{t}$ number of $K_t$'s. Since $G'$ is not a counter-example, $G'$ has exactly the optimal number of $K_t$'s and thus is one of the extremal graphs. However, adding the clique back implies that $G$ is also one of the extremal graphs, a contradiction. 
	
Now, suppose that $G$ contains a proper subgraph $H$ that is a union of connected components of $G$ where $|V(H)| \ge k-1$. Since $H$ is not a counter-example to Theorem~\ref{path}, either $H$ contains strictly less $K_t$'s than an extremal graph with the same number of vertices, or $H$ is an extremal graph. In the first case, replacing $H$ with one of the extremal graphs results in a graph that strictly increases the number of $K_t$'s in $G$ while maintaining the same number of vertices, creating a counter-example with the same number of vertices as $G$ but with more copies of $K_t$ than $G$, a contradiction. In the latter case, $H$ contains at least one copy of $K_{k-1}$, contradicting our previous claim. Thus, we may assume that all proper subgraphs that are a union of connected components of $G$ have strictly less than $k-1$ vertices. 
	
Let $G_1$ be a connected component of $G$ and $G_2=G\backslash G_1$. For $i =1,2$, let $r_i$ denote the number of vertices in $G_i$. By the previous claim, we know that both $r_1$ and $r_2$ are strictly less than $k-1$. Furthermore, as we have discussed at the beginning of this section, we can assume that the number of vertices in $G$ is at least $k$, i.e., $r_1 + r_2 \ge k$. Thus, from Lemma~\ref{easy:convex}, by setting $x = k-1$, $y = r_1$, $z = r_2$ and $w = r_1 + r_2 - (k-1)$, we obtain the following: $$\binom{k-1}{t} + \binom{r_1+r_2-(k-1)}{t} > \binom{r_1}{t} + \binom{r_2}{t} \ge k_t(G).$$ Thus, $K_{k-1} \cup K_{r_1+r_2-(k-1)}$ is a $P_k$-free graph with the same number of vertices as $G$, but with strictly more copies of $K_t$, a contradiction. 
\end{proof}

\begin{lemma} \label{mindeg}
If $G$ is a minimum counter-example, then $G$ has minimum degree at least $\max(r, t-1)$.
\end{lemma}

\begin{proof}
If $G$ is a counter-example containing a vertex $v$ with degree less than $t-1$, then, one can remove $v$ from $G$ to obtain a smaller counter-example $G \backslash v$, which contradicts our choice of $G$. Thus, the minimum degree of $G$ is at least $t-1$.
	
To show that the minimum degree of $G$ is at least $r$, for the sake of contradiction, we assume that there exists a vertex $v$ with degree $\b < r$. We may trivially assume that $r>0$. Note that there are at most $\binom{\b}{t-1}$ copies of $K_t$ containing $v$. Consider $G\backslash v$, which is not a counter-example to Theorem~\ref{path}. Since $\b < r$,
	\begin{align*}
		k_t(G\backslash v) &\le q\binom{k-1}{t}+\binom{r-1}{t}\\
		&= q\binom{k-1}{t}+\binom{r}{t}-\binom{r-1}{t-1}\\
		&\le k_t(G)-\binom{\b}{t-1} \le k_t(G\backslash v)	
	\end{align*}

Thus, the inequalities above are tight and $G\backslash v$ is one of the extremal graphs. In particular, $G\backslash v$ and $G$ contain $q$ copies of $K_{k-1}$. Note that all these cliques are disjoint and disconnected; otherwise, $G$ contains a copy of $P_k$. Since $n\ge k, q\ge 1$ and thus $G$ is not connected, contradicting Lemma~\ref{con}. 

\end{proof}

Unlike the minimum degree, we are unable to bound the maximum degree of $G$ in such a way. However, we can obtain an upper bound for the number of $K_t$'s containing any vertex $v$ based on its degree and $k$. Slightly abusing the notation, let $k_t(v)$ denote the number of $K_t$'s containing $v$. 

\begin{lemma} \label{smallcase}
Let $G$ be a minimum counter-example and $v$ be a vertex in $G$. Let $a$ and $b < k-2$ be non-negative integers such that $d(v) = a (k-2) + b$. Then, $k_t(v) \le a \binom{k-2}{t-1} + \binom{b}{t-1}$. Moreover, if $d(v) \ge k-1$, then the inequality is strict, i.e., equality cannot be achieved.
\end{lemma}

\begin{proof}
Let $N(v)$ denote the set of all neighbors of $v$. Note that the number of copies of $K_t$ containing $v$ is the same as the number of copies of $K_{t-1}$ in $N(v)$. Since $G$ is $P_k$-free, it follows that the graph induced by $N(v)$ is $P_{k-1}$-free. Hence, by the minimality of $t$, we may apply Theorem~\ref{path} on $N(v)$ for $t \ge 4$ or apply Theorem~\ref{erdos} on $N(v)$ for $t=3$, and obtain $k_t(v) \le a \binom{k-2}{t-1} + \binom{b}{t-1}$. Thus, it remains to show that this bound is never achieved when $d(v) \ge k-1$. Suppose for the sake of contradiction that this bound is achieved, then the graph induced by $N(v)$ is one of the extremal graphs. Since $d(v)\ge k-1$, $a\ge 1$ and the graph induced by $N(v)$ contains a copy of $P_{k-2}$. Since $d(v)\ge k-1$, there exists a vertex $u \in N(v)$ outside of this path. Then, one can extend the path to $u$ via $v$, obtaining a path on $k$ vertices, a contradiction.
\end{proof}

In the extremal graphs of Theorem~\ref{path}, the number of copies of $K_t$ is the same as $\frac{1}{t}\sum_{v}\binom{d(v)}{t-1}$. Then in order to compare $k_t(G)$ to the optimal number, it is useful to similarly construct a sequence $\{y_1, \ldots, y_N\}$ such that $k_t(G)\le \frac{1}{t}\sum_{i=1}^N \binom{y_i}{t-1}$. Motivated by this fact and the result of Lemma~\ref{smallcase}, we introduce the following definition:

\begin{definition} \label{split}
Given non-negative integers $\Delta$ and $d$, let $\spl_\D(d)$ output a sequence of $a_d$ integers of value $\D$ and one integer of value $b_d$ where $a_d\D+b_d=d$ and $b_d<\D$. Given a sequence of non-negative integers $\bar{d}=\{d_1, \ldots, d_n\}$, let $\spl_\D{\bar{d}}$ be the concatenation of the $\spl_\D(d_i)$'s, arranged in descending order.
\end{definition}

Then, we obtain the following corollary from Lemma~\ref{smallcase}:

\begin{corollary} \label{countkt}
Let $G$ be a minimum counter-example. Let $\{y_i\}_{i=1}^{N}=\spl_{k-2}(\bar{d})$ where $\bar{d}$ is the degree sequence of $G$. Then $k_t(G)\le \frac{1}{t}\sum_{i=1}^N\binom{y_i}{t-1}$. 
\end{corollary}

Now, we can compare $\sum_{i=1}^N \binom{y_i}{t-1}$ with $q(k-1) \binom{k-2}{t-1} + r \binom{r-1}{t-1} = t \cdot k_t\left(qK_{k-1}\cup K_r\right)$ using the well-known Karamata's inequality:

\begin{lemma}[Karamata's inequality \cite{K}] \label{karma}
Let $f$ be a real-valued convex function defined on $\N$. If $x_1, x_2, \ldots, x_N \in \N$ and $y_1, y_2, \ldots, y_N \in \N$ are such that 
\begin{align}
&\bullet \;\; x_1 \ge\cdots \ge x_N \text{	and	} y_1 \ge y_2 \ge \cdots \ge y_N, \label{K1}
\\ &\bullet \;\; x_1 + x_2 + \cdots x_i \ge y_1 + y_2 + \cdots + y_i \text{	for all	} i \in \{1, 2, \ldots, N-1\}, \text{	and	} \label{K2}
\\ &\bullet \;\; x_1 + x_2 + \cdots x_N = y_1 + y_2 + \cdots + y_N, \label{K3}
\end{align}
then 
\begin{equation} \label{karamata}
f(x_1) + f(x_2) + \cdots + f(x_N) \ge f(y_1) + f(y_2) + \cdots + f(y_N).
\end{equation}
\end{lemma}

\begin{definition}
For $x_1, x_2, \ldots, x_N \in \N$ and $y_1, y_2, \ldots, y_N \in \N$, we say that $(x_1, x_2, \ldots, x_N)$ majorizes $(y_1, y_2, \ldots, y_N)$ if the conditions in \eqref{K1}, \eqref{K2}, and \eqref{K3} are satisfied.
\end{definition}

\begin{proof} [Proof of Theorem~\ref{path}]
	
Let $\{y_i\}_{i=1}^N=\spl_{k-2}(\bar{d})$, where $\bar{d}$ is the degree sequence of $G$. Then, construct the following sequence, $\{x_i\}_{i=1}^N$, of the same length $N$:

\begin{itemize}
\item $x_i = k-2$ for all $1 \le i \le q(k-1)$, 
\item $x_i = r-1$ for all $q(k-1) < i \le q(k-1) + r$, and 
\item $x_i = 0$ for all $i > q(k-1) + r$.
\end{itemize}

Next we show that $(x_1, x_2, \ldots, x_N)$ majorizes $(y_1, y_2, \ldots, y_N)$. Condition~\eqref{K1} is trivially true. By applying Theorem~\ref{erdos}, the number of edges in $G$ is at most $q\binom{k-1}{2} + \binom{r}{2}$. This implies $\sum_{i=1}^N y_i\le \sum_{i=1}^N x_i$. We modify the sequences slightly to satisfy Condition~\eqref{K3}. We consider the sequence $\{y'_i\}_{i=1}^{N'}$ by adding as many $1$'s as needed at the end of $\{y_i\}_{i=1}^N$, until the condition is satisfied. We similarly create $\{x'\}_{i=1}^{N'}$ by padding the $x$ sequence with $0$'s until the two sequences have the same length. Note that Condition~\eqref{K1} remains satisfied. In order to show that Condition~\eqref{K2}, is satisfied, we consider the following three ranges separately.
\smallskip

Case 1: If $1 \le i \le q(k-1)$, then $$x'_1 + x'_2 + \cdots + x'_i = i (k-2) \ge y'_1 + y'_2 + \cdots + y'_i.$$
\smallskip

Case 2: $q(k-1) + r < i < N'$, then $$x'_1 + x'_2 + \cdots + x'_i = x'_1 + x'_2 + \cdots + x'_{N'} = y'_1 + y'_2 + \cdots + y'_{N'} \ge y'_1 + y'_2 + \cdots + y'_i.$$
\smallskip

Case 3: If $q(k-1) < i \le q(k-1) + r$, then since every vertex in $G$ has degree at least $r$, $\spl_{k-2}(d(v))$ outputs an integer at least $r$ for all $v \in V(G)$. Then, it follows that the first $n=q(k-1) + r$ elements in the sequence $y'$ are all at least $r$. Hence, $\sum_{j=i+1}^{q(k-1)+r} x'_j \le \sum_{j=i+1}^{q(k-1)+r} y'_j$. Thus,
\begin{align*}
&\sum_{j=1}^{N'} x'_i - \sum_{j=i+1}^{q(k-1)+r} x'_j \ge \sum_{j=1}^{N'} y'_i - \sum_{j=i+1}^{q(k-1)+r} y'_j \\
\Longrightarrow \;\; &\sum_{j=1}^{q(k-1)+r} x'_i - \sum_{j=i+1}^{q(k-1)+r} x'_j \ge \sum_{j=1}^{q(k-1)+r} y'_i - \sum_{j=i+1}^{q(k-1)+r} y'_j \\
\Longrightarrow \;\; &x'_1 + x'_2 + \cdots + x'_i \ge y'_1 + y'_2 + \cdots + y'_i.
\end{align*}
 
This proves $(x'_1, \ldots, x'_{N'})$ majorizes $(y'_1, \ldots, y'_{N'})$. Then, apply Karamata's inequality~\eqref{karamata} for the convex function $f : \N \rightarrow \N$ defined by $f(x) = \binom{x}{t-1}$ to obtain 
\begin{equation} \label{contradictionkt}
\sum_{i=1}^N \binom{y_i}{t-1} = \sum_{i=1}^{N'} \binom{y_i'}{t-1} \le \sum_{i =1}^{N'}\binom{x_i'}{t-1} = \sum_{i=1}^N\binom{x_i}{t-1}=t \cdot k_t(qK_{k-1}\cup K_r).
\end{equation}

If $G$ has a vertex of degree greater than $k-2$, then it follows from the second part of Lemma~\ref{smallcase} and Corollary~\ref{countkt} that $\frac{1}{t}\sum_{i=1}^N \binom{y_i}{t-1}> k_t(G)\ge k_t(qK_{k-1}\cup K_r)$, a contradiction to \eqref{contradictionkt}. Similarly, when the maximum degree of $G$ is at most $k-2$, the inequalities in \eqref{contradictionkt} and Corollary~\ref{countkt} should be equalities. Consequently, for every vertex $v$, we have that $k_t(v) = \binom{d(v)}{t-1}$. Using Lemma~\ref{mindeg}, we know that every vertex $v$ satisfies that $d(v) \ge t-1$. Hence, it follows that for every vertex $v$, the neighborhood $N(v)$ induces a clique. This fact together with Lemma~\ref{con} implies that $G$ must be a complete graph. Since $n\ge k$, $G$ contains a copy of $P_k$, a contradiction. 

\end{proof}

\section{Proof of Theorem~\ref{edge}}
\label{sec:newedge}

For fixed $t\ge 3$, we assume that $G$ is a minimum counter-example to Theorem~\ref{edge}. In particular, we assume that $G$ has $q\binom{k-1}{2}+\binom{r}{2}+s$ edges where $0\le s< r\le k-2$, $k_t(G)\ge q\binom{k-1}{t}+\binom{r}{t}+\binom{s}{t-1}$, and $G$ is not one of the extremal graphs (i.e. $G \not\in \Gf(qK_{k-1}, L)$, where $L\in \C_{t, k}(\binom{r}{2}+s)$). We may assume that $0\le s<r$ because any non-negative integer $b$ can be uniquely written in the form $b = \binom{r}{2} + s$ for integers $0 \le s < r$. Furthermore, we assume that the number of copies of $K_t$ in any other graph with the same number of edges is at most $k_t(G)$. Note that if $q=0$, then Theorem~\ref{edge} follows from Theorem~\ref{maxkt}. Thus, we may also assume that $q \ge 1$.

\begin{lemma} \label{bigclique}
If $G$ is a minimum counter-example, then $G$ does not contain a copy of $K_{k-1}$.
\end{lemma}

\begin{proof}
Suppose for the sake of contradiction that $G$ does contain a copy of $K_{k-1}$. Fix one such copy $C$ of $K_{k-1}$, let $v_1, v_2, \ldots, v_{k-1}$ denote the vertices of the clique $C$. We first claim that for all $i \neq j$, there is no path $P$ with at least three vertices between the vertices $v_i$ and $v_j$ such that all internal vertices are outside of $C$. Indeed, if this is not the case, then we can easily find a cycle of length at least $k$ by using the path $P$ and the clique $C$, contradicting the fact that $G$ is $C_{\ge k}$-free. Then, it follows that for any edge $e\in E(C)$, all the copies of $K_t$ using $e$ must be subgraphs of $C$. Then, by removing the edges in $C$, we obtain a graph $G'$ with $(q-1)\binom{k-1}{2}+\binom{r}{2}+s$ edges and contains at least $(q-1)\binom{k-1}{t}+\binom{r}{t}+\binom{s}{t-1}$ number of $K_t$'s. Since $G'$ is not a counter-example, $G'$ has exactly the optimal number of $K_t$'s and thus is one of the extremal structures where $G'\in \mathcal{G}_{forest}((q-1)K_{k-1}, L)$ for some $L\in \C_{t, k}(\binom{r}{2}+s)$. 

To achieve a contradiction, we show that $G\in \mathcal{G}_{forest}(qK_{k-1}, L)$. Let $F_1, F_2, \ldots, F_{q-1}$ be the $q-1$ edge-disjoint copies of $K_{k-1}$ in $G'$. Let $L_1, L_2, \ldots, L_{q'}$ be the connected components in $L$. Let $\mathcal{F}=\{F_1, \dots, F_{q-1}, L_1, \ldots, L_{q'}\}$. Note that the edge sets of the graphs in $\mathcal{F}$ and $E(C)$ partition the edges of $G$. Furthermore, $|V(G_1)\cap V(G_2)|\le 1$ for all distinct $G_1, G_2\in \mathcal{F}$. For any $G_1\in \mathcal{F}$, $|V(G_1)\cap V(C)|\le1$; otherwise, two of the vertices in $C$ has a path that does not use any edges in $C$, contradicting our previous claim about $C$. Then, $G$ satisfies the first condition in Definition~\ref{def:family} with respect to $\mathcal{G}_{forest}(qK_{k-1}, L)$. Let $H$ be the auxiliary graph created from $G$ according to the second condition in Definition~\ref{def:family}. Since $G'\in \mathcal{G}_{forest}(\mathcal{F})$ and no two vertices in $C$ have an external path between them outside of $C$ in $G$, it follows that $H$ is a forest. Then, $G\in \mathcal{G}_{forest}(\mathcal{F}, C)$ and thus $G$ is one of the extremal graphs, a contradiction. 
\end{proof}

\begin{lemma} \label{connected}
If $G$ is a minimum counter-example, then $G$ is 2-connected.
\end{lemma}

Lemma~\ref{connected} is proved at the end of this section to avoid being sidetracked by simple but long convexity arguments. On an intuitive level, if a counter-example $G$ is not connected, then one of its connected components should be a smaller counter-example. Similarly, if a counter-example $G$ is connected but not 2-connected, then $G$ has a cut vertex $v$. Then, by splitting the vertex $v$, we should be able to find a smaller minimal counter-example. For now, we take this fact for granted and prove other structural properties about $G$.

\begin{lemma}
\label{mindegedge}
If $G$ is a minimum counter-example, then $G$ has minimum degree at least $r$.
\end{lemma}

\begin{proof}
For the sake of contradiction, assume that there exists $v\in V(G)$ whose degree is $0 < \b < r$. Consider the graph $G\backslash v$. Let $k=k_t(G)-k_t(G\backslash v)$ be the number of $K_t$'s that contain $v$. Note that $k \le \binom{\b}{t-1}$. 
	
Since $G\backslash v$ is not a counter-example, $k_t(G\backslash v)\le q\binom{k-1}{t}+k_t(L_{\binom{r}{2}+s-\b})$. If $\b\le s$, then by Lemma~\ref{easy:convex}, $k_t(L_{\binom{r}{2}+s}) - k_t(L_{\binom{r}{2}+s-\b}) = \left(\binom{r}{t}+\binom{s}{t-1}\right)-\left(\binom{r}{t} + \binom{s-\b}{t-1}\right)\ge \binom{\b}{t-1}$. If $\b>s$, since $\b\le r-1$, then $k_t(L_{\binom{r}{2}+s}) - k_t(L_{\binom{r}{2}+s-\b}) = \left(\binom{r}{t}+\binom{s}{t-1}\right)-\left(\binom{r-1}{t} + \binom{r-1+s-\b}{t-1}\right)\ge \binom{\b}{t-1}$. Note that in both cases, $$k_t(L_{\binom{r}{2}+s}) - k_t(L_{\binom{r}{2}+s-\b}) \ge \binom{\b}{t-1}.$$ 
	
Then, we have the following.
	
	\begin{align*}
	k_t(G\backslash v) &= k_t(G) - k \\
	&\ge q\binom{k-1}{t} + k_t(L_{\binom{r}{2}+s}) - \binom{\b}{t-1} \\
	&\ge q\binom{k-1}{t} + k_t(L_{\binom{r}{2}+s-\b}) \ge k_t(G\backslash v).
	\end{align*}
	
Since $G\backslash v$ is not a counter-example, we conclude that it is one of the extremal graphs in Theorem~\ref{edge}. Thus, since $q\ge 1$, both $G\backslash v$ and $G$ contain a copy of $K_{k-1}$, contradicting Lemma~\ref{bigclique}. 
	
\end{proof}

\begin{lemma} \label{smallr}
If $G$ is a minimum counter-example, then $r \le \frac{k-1}{2}$.
\end{lemma}

To prove this lemma, we require the following well-known theorem by Dirac \cite{D}.

\begin{theorem} [\cite{D}] \label{dirac}
If $G$ is an $n$-vertex 2-connected graph with minimum-degree at least $r$, then there exists a cycle of length at least $\min(n, 2r)$. 
\end{theorem}

\begin{proof} [Proof of Lemma~\ref{smallr}]
This follows immediately from Lemmas~\ref{connected} and~\ref{mindegedge}, and Theorem~\ref{dirac}. 
\end{proof}

Next, we bound the number of vertices in $G$. 

\begin{lemma} \label{lem:largen}
If $G$ is a minimum counter-example, then $G$ has at least $q(k-1) + (r+1)$ vertices.
\end{lemma}

To prove this lemma, we make use of a classical result by Whitney \cite{W32} and a recent result by Luo \cite{Luo}. 

\begin{theorem} [\cite{W32}] \label{Whitney}
A connected graph $G$ with at least three vertices is 2-connected if and only if for every two vertices $u, v \in V(G)$, there is a cycle containing both $u$ and $v$.
\end{theorem}

\begin{theorem} [\cite{Luo}] \label{Luo3}
Let $n \ge k \ge 5$. If $G$ is a 2-connected $n$-vertex graph with no cycle of length at least $k$, then $k_t(G) \le \max\{f_t(n,k,2), f_t(n,k,\lfloor\frac{k-1}{2}\rfloor)\}$, where $f_t(n,k,a) = \binom{k-a}{t} + (n-k+a)\binom{a}{t-1}$.
\end{theorem}

\begin{proof} [Proof of Lemma~\ref{lem:largen}]
First, we claim that $k\ge 5$. Suppose for the sake of contradiction that $k\le 4$. By Lemma~\ref{connected} and Theorem~\ref{Whitney}, any two vertices $u, v$ are in a cycle. Since $G$ is $C_{\ge k}$-free, then every pair of vertices $u, v$ is in a triangle and in particular, edge $uv$ exists. Therefore $G$ is a clique. However, since $G$ is $C_4$-free, $G$ has at most three vertices and is trivially one of the extremal graphs, a contradiction.
	
Suppose for the sake of contradiction that $G$ has at most $n=q(k-1) + r$ vertices. In order to get a contradiction, it is enough to prove that $k_t(G) < T_t:= k_t(qK_{k-1} \cup L_{\binom{r}{2} + s})=q\binom{k-1}{t} + \binom{r}{t} + \binom{s}{t-1}$. By Lemma~\ref{connected} and Theorem~\ref{Luo3}, it suffices to prove that $\max\{f_t(n,k,2), f_t(n,k,\lfloor\frac{k-1}{2}\rfloor)\} < T_t$.

We first show that $f_t(n,k,2) < T_t$ as follows:
\begin{align*}
f_t(n,k,2) &= \binom{k-2}{t} + ((q-1)(k-1)+r+1)\binom{2}{t-1} \\
&\le \binom{k-2}{t} + ((q-1)(k-1)+r+1) \\
&< q\binom{k-1}{t} + \binom{r}{t} + \binom{s}{t-1} = T_t.
\end{align*}

Now, we show that $f_t(n,k,\lfloor\frac{k-1}{2}\rfloor) = \binom{k-\lfloor\frac{k-1}{2}\rfloor}{t} + (q(k-1)+r-k+\lfloor\frac{k-1}{2}\rfloor)\binom{\lfloor\frac{k-1}{2}\rfloor}{t-1} < T_t$. Note that it suffices to prove the following:
$$\binom{\lceil\frac{k+1}{2}\rceil}{t} + \left((q-1)(k-1)+r+\left\lfloor\frac{k-3}{2}\right\rfloor\right)\binom{\lfloor\frac{k-1}{2}\rfloor}{t-1} < q\binom{k-1}{t} .$$

By rewriting $q\binom{k-1}{t}$ as $(q-1)\binom{k-1}{t}+\binom{k-1}{t}$ and using the fact that $r\le \lfloor \frac{k-1}{2}\rfloor$ (by Lemma~\ref{mindegedge}), the above inequality is implied by the following two inequalities:

\begin{align}
(q-1)(k-1)\binom{\lfloor\frac{k-1}{2}\rfloor}{t-1} &\le (q-1)\binom{k-1}{t} \label{ineq1} \\
\binom{\lceil\frac{k+1}{2}\rceil}{t} + 2\left\lfloor\frac{k-1}{2}\right\rfloor\binom{\lfloor\frac{k-1}{2}\rfloor}{t-1} & < \binom{k-1}{t} \label{ineq2}
\end{align}

To prove \eqref{ineq1}, we may assume that $\lfloor\frac{k-1}{2}\rfloor\ge t-1$, otherwise the inequality trivially holds. Rearranging the terms, it is equivalent to show $\prod_{i=0}^{t-2} \frac{k-2-i}{\lfloor\frac{k-1}{2}\rfloor -i} \ge t$. For the left-hand side, when $i=0, \frac{k-2}{\lfloor \frac{k-1}{2}\rfloor} \ge 1.5$ since $k\ge 5$. For the subsequent fractions, note that for $1\le i\le t-2$,  $\frac{k-2-i}{\lfloor\frac{k-1}{2}\rfloor-i} \ge \frac{k-3}{\lfloor\frac{k-1}{2}\rfloor -1}$. Since $k\ge 5$, these fractions are at least $2$. Then, it suffices to show that $1.5(2)^{t-2}\ge t$. This holds for $t\ge 3$, proving the inequality of \eqref{ineq1}.

To prove \eqref{ineq2}, we require the following inequality:
$$\binom{x}{t} + y\binom{x}{t-1} + x\binom{y}{t-1} \le \binom{x+y}{t}.$$

If the above inequality holds, \eqref{ineq2} follows immediately by setting $x=\lceil\frac{k+1}{2}\rceil$ and $y=\lfloor\frac{k-1}{2}\rfloor$. To prove the inequality, consider the following combinatorial argument. 

Suppose there are two groups of people $A$ and $B$ of size $x$ and $y$, respectively. The right-hand side counts the number of ways to form a team of $t$ people out of $A\cup B$. Possible ways of forming a $t$-member team include: choosing all $t$ members from $A$, choosing one member from $A$ while the rest from $B$, and choosing one member from $B$ while the rest from $A$. Thus, the left-hand side is at most the right-hand side. 
\end{proof}

We require one more lemma before proving our main theorem. The following allows us to better compare the value $k_t(G)$ and the optimal value $T_t$. 

\begin{lemma} \label{cor}
If $G$ is a minimum counter-example with $m$ edges, then $$k_t(G) \le \frac{1}{t} \left[\a \binom{k-2}{t-1} + (n-\a) \binom{r}{t-1}\right],$$ where $\a$ is a real number satisfying $\a(k-2) + (n-\a)r = 2m$ and $n = q(k-1) + (r+1)$.
\end{lemma}

\begin{proof}
Let $G$ be a minimum counter-example. Keep in mind that the number of vertices in $G$ is at least $n$ by Lemma~\ref{lem:largen}. For any vertex $v$ in a $C_k$-free graph, the neighborhood $N(v)$ is $P_{k-1}$-free, and hence by Theorem~\ref{path} (or Theorem~\ref{erdos} for $t=3$), $k_t(v) \le a \binom{k-2}{t-1} + \binom{b}{t-1}$, where $a$ and $b$ are non-negative integers satisfying $d(v) = a(k-2) + b$ and $b<k-2$. By using Definition~\ref{split}, we have that $\sum_{v \in V(G)} k_t(v) \le \sum_{i=1}^N \binom{y_i}{t-1}$, where $\{y_i\}_{i=1}^{N}=\spl_{k-2}(\bar{d})$ and $\bar{d}$ is the degree sequence of $G$. Now, define the sequence $\{x_i\}_{i=1}^n$ such that $\sum_{i=1}^n x_i = 2m$ and there exists $j$ such that $x_i=k-2$ for all $1\le i\le j$, $r \le x_{j+1} \le k-2$ and $x_i=r$ for all $j+1<i\le n$. Similar to the proof of Theorem~\ref{path} in Section 2, Karamata's inequality~\ref{karma} implies that $\sum_{v \in V(G)} k_t(v) \le \sum_{i=1}^n \binom{x_i}{t-1}$. Now, find the unique $0 \le \a' \le 1$ such that $\a' (k-2) + (1-\a') r = x_{j+1}$, and apply Jensen's inequality to get $\binom{x_{j+1}}{t-1} \le \a' \binom{k-2}{t-1} + (1-\a') \binom{r}{t-1}$. Our lemma follows immediately with $\a=j+\a'$.
\end{proof}

Now, we are ready to prove Theorem~\ref{edge}. 

\begin{proof}[Proof of Theorem~\ref{edge}]
Let $G$ be a minimum counter-example. By Lemma~\ref{cor}, one can check that $\a = q(k-1) - \frac{2(r-s)}{k-2-r}$ is the desirable $\alpha$. Then, we have the following:
	\begin{align*}
	k_t(G) &\le \frac{1}{t} \left[ \left(q(k-1) - \frac{2(r-s)}{k-2-r}\right) \binom{k-2}{t-1} + \left(r+1+\frac{2(r-s)}{k-2-r}\right) \binom{r}{t-1} \right] \\
	&\le T_t -\frac{2(r-s)}{t(k-2-r)} \left(\binom{k-2}{t-1} - \binom{r}{t-1}\right) + \binom{r}{t-1} - \binom{s}{t-1}. 
	\end{align*}
	
Let $C_1=\frac{2(r-s)}{t(k-2-r)} \left(\binom{k-2}{t-1} - \binom{r}{t-1}\right)$ and $C_2=\binom{r}{t-1} - \binom{s}{t-1}$. In order to show that $k_t(G) < T_t$, it suffices to prove that $C_2-C_1<0$. Since $r>s$, the inequality is strict when $r<t-1$. Thus, we may assume $r\ge t-1$. In order to finish the proof, we claim that it suffices to prove the following two inequalities:
	
	\begin{equation} \label{eq1}
		\binom{r}{t-1} - \binom{s}{t-1} \le (r-s)\binom{r-1}{t-2}	
	\end{equation}
	
		\begin{equation} \label{t}
	\frac{\binom{k-2}{t-1} - \binom{r}{t-1}}{(k-2-r) \binom{r-1}{t-2}} > \frac{t}{2}.
	\end{equation}
	
If the inequalities above were true, Equation~\eqref{t} provides a bound to $C_1$. Then, $C_2-C_1<0$ follows from Equation~\eqref{eq1}.
	
To show that Equation~\eqref{eq1} holds, the left-hand side can be interpreted as choosing $t-1$ people in a group $R$ of size $r$ without choosing all of them from a subgroup $S\subset R$ of size $s$. The right-hand side is an upper bound on the number of ways a team can be formed by making sure at least one of them is from $R \setminus S$. It follows that the right-hand side is at least as large as the left-hand side.

To prove Equation~\eqref{t}, since $t\ge 3$, we can first lower-bound the numerator of the left-hand side by the following:
	\begin{equation} \label{4}
	\binom{k-2}{t-1} - \binom{r}{t-1} \ge (k-2-r) \binom{r}{t-2} + \binom{k-2-r}{2} \binom{r}{t-3}.
	\end{equation}
	
The above inequality can be proved by using a similar argument as before: the left-hand side is forming a team of $t-1$ players from a group $K$ of size $k-2$ without having all of them from a subgroup $R\subset K$ of size $r$. Meanwhile, the right-hand side corresponds to first choosing one or two people from $K\setminus R$ and filling the rest of the team from $R$. Then,

	\begin{align} \label{6}
	\frac{\binom{k-2}{t-1} - \binom{r}{t-1}}{(k-2-r) \binom{r-1}{t-2}} &\ge \frac{(k-2-r)\binom{r}{t-2}+\binom{k-2-r}{2}\binom{r}{t-3}}{(k-2-r)\binom{r-1}{t-2}} \nonumber \\
	&= \frac{r}{r-t+2} + \frac{r}{r-t+3}\cdot \frac{k-3-r}{r-t+2} \cdot \frac{t-2}{2}.
	\end{align}
	
Thus, for $t=3$, one can easily plug in the value of $t$ in \eqref{6} and check the validity of \eqref{t} by using Lemma~\ref{smallr}. For $t \ge 4$, we have the following from \eqref{6}:
	
	 \begin{align*} 
	\frac{\binom{k-2}{t-1} - \binom{r}{t-1}}{(k-2-r) \binom{r-1}{t-2}} > 1 + \frac{k-3-r}{r-2} \cdot \frac{t-2}{2}. 
	\end{align*}
	
Thus, to prove \eqref{t}, it is enough to show that $\frac{k-3-r}{r-2} \ge 1$. This follows from Lemma~\ref{smallr}, which concludes our proof of Theorem~\ref{edge}.
	
\end{proof}

We now finish this section by giving a proof of Lemma~\ref{connected}. We need the following couple of lemmas about a minimum counter-example to Theorem~\ref{edge}. 

\begin{lemma}
\label{lem:posr}
If $G$ is a minimum counter-example, then $r \ge 2$. 
\end{lemma}
	
\begin{proof}
Suppose for the sake of contradiction that $r=1$ and $s=0$. This implies that $|E(G)|=q\binom{k-1}{2}$ and $k_t(G)\ge q\binom{k-1}{t}$. We split into two cases depending on the maximum degree of $G$. 
\smallskip
	
Case 1: If the maximum degree of $G$ is strictly less than $k-2$, then for any edge $e$, its endpoints have less than $k-3$ common neighbors in $G$. Hence, $e$ belongs in less than $\binom{k-3}{t-2}$ distinct copies of $K_t$. Summing over all edges, since each $K_t$ is counted $\binom{t}{2}$ times, $k_t(G)$ is less than $q\binom{k-1}{2}\frac{\binom{k-3}{t-2}}{\binom{t}{2}}=q\binom{k-1}{t}$, a contradiction. 
\smallskip
	
Case 2: If the maximum degree of $G$ is at least $k-2$, we claim that $G$ contains a copy of $K_{k-1}$ which contradicts Lemma~\ref{bigclique}. For a vertex $v$, by Theorem~\ref{path} and using the fact that the neighborhood of $v$ is $P_{k-1}$-free, we have that 
	\begin{equation} \label{copy}
	k_t(v) \le a \binom{k-2}{t-1} + \binom{b}{t-1},
	\end{equation}
	where $a$ and $b < k-2$ are non-negative integers satisfying $d(v) = a(k-2) + b$. By using the convexity of the function $x \rightarrow \binom{x}{t-1}$ with the set of integers as domain, it is easy to see that $a\binom{k-2}{t-1} + \binom{b}{t-1} \le \frac{d(v)}{k-2} \binom{k-2}{t-1}$. Hence, it follows that 
	\begin{align}
	k_t(G) = \frac{1}{t} \sum_{v \in V(G)} k_t(v) 
	\le \frac{1}{t} \sum_{v \in V(G)} \frac{d(v)}{k-2} \binom{k-2}{t-1} 
	= q \binom{k-1}{t}. \label{conseq}
	\end{align} 
Thus, we have an equality in \eqref{conseq} and consequently, we have equality in \eqref{copy} for each vertex of $G$. Suppose that $v$ is a vertex with degree at least $k-2$. Then, by using the equality in \eqref{copy} for $v$, it follows from Theorem~\ref{path} that there is a copy of $K_{k-2}$ in the neighborhood $N(v)$ of $v$. This $K_{k-2}$ together with the vertex $v$ makes a copy of $K_{k-1}$ in $G$, contradicting Lemma~\ref{bigclique}.
\end{proof}

\begin{lemma}
\label{alledget}
If $G$ is a minimum counter-example, then every edge in $G$ is in a copy of $K_t$. 
\end{lemma}

\begin{proof}
Suppose for the sake of contradiction that an edge $e$ is not in a copy of $K_t$. Consider the graph $G'=G\backslash e$. 
\smallskip
	
Case 1: $s>0$. Observe that $q\binom{k-1}{t}+\binom{r}{t}+\binom{s}{t-1}\le k_t(G) = k_t(G') \le q\binom{k-1}{t}+\binom{r}{t}+\binom{s-1}{t-1}$. Since the right-hand side is at most the left-hand side, the above equation is at equality and $s<t-1$. Then, $G'$ has an extremal structure $\Gf(qK_{k-1}, L)$, where $L\in \C_{t, k}(\binom{r}{2}+s-1)$. Since $q\ge 1$, $G$ contains a copy of $K_{k-1}$, contradicting Lemma~\ref{bigclique}.
\smallskip
	
Case 2: $s=0$. Keep in mind that by Lemma~\ref{lem:posr}, we have that $r-2\ge 0$. Then, $q\binom{k-1}{t}+\binom{r}{t}\le k_t(G) = k_t(G') \le q\binom{k-1}{t}+\binom{r-1}{t}+\binom{r-2}{t-1}$. Once again, the above equation is tight and thus $r<t$. Then, $G'$ is an extremal graph with $q$ cliques of order $k-1$ and a graph $L\in \C_{t, k}(\binom{r}{2}-1)$. Similarly, $G$ contains a copy of $K_{k-1}$, contradicting Lemma~\ref{bigclique}.  
	
\end{proof}

\begin{proof}[Proof of Lemma~\ref{connected}]
Suppose for the sake of contradiction that $G$ is not 2-connected. Let $v$ be a vertex in $G$ such that $G' = G \backslash v$ is not connected. Now, suppose $G'$ contains a proper subgraph $H'$ that is a union of connected components of $G'$ such that $|E(H)| \ge \binom{k-1}{2}$ where $H$ is the graph induced by $V(H')\cup\{v\}$. Since $H$ is not a counter-example to Theorem~\ref{edge}, either $H$ contains strictly less $K_t$'s than an extremal graph with the same number of edges, or $H$ is an extremal graph. In the first case, deleting the edges in $H$ and replacing with one of the extremal graphs results in a graph that strictly increases the number of $K_t$'s in $G$ while maintaining the same number of edges, creating a counter-example with the same number of edges as $G$ but with more copies of $K_t$, a contradiction. In the latter case, $H$ contains at least one copy of $K_{k-1}$, contradicting Lemma~\ref{bigclique}. Thus, we may assume that for all proper subgraphs $H'$ that is a union of connected components of $G' = G \backslash v$, the graph induced by $V(H') \cup \{v\}$ has strictly less than $\binom{k-1}{2}$ edges. 
	
Let $G'_1$ be a connected component of $G'$ and $G'_2=G'\backslash G'_1$. Let $G_i$ be the graph induced by $V(G'_i) \cup \{v\}$ for $i=1, 2$. Note that if one of $G_1$ and $G_2$ is not one of the extremal structures, then by replacing them with an extremal structure, one strictly increases the number of $K_t$'s, contradicting the choice of $G$. Since neither are counter-examples nor contain at least $\binom{k-1}{2}$ edges, we may assume that $G_i\in \C_{t, k}(\binom{r_i}{2}+s_i)$, where $|E(G_i)|=\binom{r_i}{2}+s_i$ and $0\le s_i< r_i\le k-2$ for $i=1, 2$. Observe that if $0<s_i<t-1$, some edges in $G_i$ will not be part of any $K_t$'s, contradicting Lemma~\ref{alledget}. A similar contradiction is achieved if $r_i<t$. Thus we may assume that $r_i\ge t$ and either $s_i=0$ or $s_i\ge t-1$ for $i=1, 2$. Observe in either case, $G_i$ is a colex graph. Without loss of generality, we can assume that $r_1 \ge r_2$. Now, depending on the values of $s_1$ and $s_2$, we move a certain amount of edges from $G_2$ to $G_1$ and obtain a graph with strictly more $K_t$'s than before, achieving a contradiction to the choice of $G$.
\smallskip
	
Case 1: $s_1, s_2> 0$. Let $s'=\min\{s_2, r_1-s_1\}$. Note that $s' \ge 1$. Let $G_1''$ be the colex graph with $|E(G_1)|+s'=\binom{r_1}{2}+s_1+s'$ edges, let $G_2''$ be the colex graph with $|E(G_2))|-s'=\binom{r_2}{2}+s_2-s'$ edges. Essentially, we are moving $s'$ edges from $G_2$ to $G_1$. The value $s'$ is chosen such that this process is equivalent to moving one edge at a time from $G_2$ to $G_1$ and stopping as soon as one of $G_1$ or $G_2$ becomes a clique. Doing so keeps the calculation of $k_t(G_1'')$ and $k_t(G_2'')$ as simple as possible. Note that $s_1, s_2 \ge t-1\ge 1$, and $s' \ge 1$. Hence, we have the conditions that $s_1, s_2 < s_1+s'$ and $s_1+s' \ge t-1$. Then, it follows from Lemma~\ref{easy:convex} that
	
	\begin{align*}
		k_t(G_1\cup G_2) &= \binom{r_1}{t}+\binom{s_1}{t-1}+\binom{r_2}{t}+\binom{s_2}{t-1}\\
		&< \binom{r_1}{t} + \binom{s_1+s'}{t-1} + \binom{r_2}{t} + \binom{s_2-s'}{t-1}\\
		&= k_t(G_1''\cup G_2''), 
	\end{align*} 
	a contradiction.
\smallskip
	
Case 2: $s_1 = 0$, $s_2> 0$. Let $s'=\min\{r_2 -1, r_1-s_2\}$. Since $r_1\ge r_2 > s_2 \ge t-1 \ge 1$, we have $s'\ge 1$. Note that after deleting $s_2$ edges from $G_2$, there are still $s'$ more edges one can remove. Then, let $G_1''$ be the colex graph with $|E(G_1)|+s_2+s'=\binom{r_1}{2}+s_2+s'$ edges and $G_2''$ be the colex graph with $|E(G_2)|-s_2-s'=\binom{r_2-1}{2} +r_2-1-s'$ edges (moving $s_2+s'$ edges from $G_2$ to $G_1$). 
Since $s_2\ge t-1$ and $s'\ge 1$, we have that $s_2 + s'\ge t-1$ and $s_2 + s' > s_2$. Furthermore, since $r_1\ge r_2$, we have that $s_2+s' = \min\{r_2-1+s_2, r_1\}> r_2-1$. Then, by Lemma~\ref{easy:convex},
	
	\begin{align*}
		k_t(G_1\cup G_2) &= \binom{r_1}{t}+\binom{r_2-1}{t}+\binom{r_2-1}{t-1}+\binom{s_2}{t-1}\\
		&< \binom{r_1}{t}  +\binom{s_2+s'}{t-1} + \binom{r_2-1}{t} + \binom{r_2-1-s'}{t-1}\\
		&= k_t(G_1''\cup G_2''),
	\end{align*}
	a contradiction. 
\smallskip
	
Case 3: $s_1>0$, $s_2 = 0$. Let $s'=\min\{r_1-s_1, r_2-1\}$. Keep in mind that $r_1>s_1\ge t-1\ge 1$. Let $G_1''$ be the colex graph with $|E(G_1)|+s'=\binom{r_1}{2}+s_1+s'$ edges and let $G_2''$ be the colex graph with $|E(G_2)|-s'=\binom{r_2-1}{2}+r_2-1-s'$ edges (moving $s'$ edges from $G_2$ to $G_1$ while ensuring not adding more than $r_1-s_1$ many edges to $G_1$). Note that $s'\ge 1$ and $s_1+s'\ge t-1$. Since $r_1\ge r_2$, we have that $s_1+s'=\min\{r_1, r_2-1+s_1\}>r_2-1$. Then, it follows from Lemma~\ref{easy:convex} that
	
	\begin{align*}
		k_t(G_1\cup G_2) &= \binom{r_1}{t}+\binom{s_1}{t-1}+\binom{r_2-1}{t}+\binom{r_2-1}{t-1}\\
		&< \binom{r_1}{t}+\binom{s_1+s'}{t-1} +\binom{r_2-1}{t}+\binom{r_2-1-s'}{t-1}\\
		&= k_t(G_1''\cup G_2''),
	\end{align*}
	a contradiction.
\smallskip
	
Case 4: $s_1=s_2=0$. Keep in mind that $r_2\ge t\ge 2$. Let $G_1''$ be the colex graph with $|E(G_1)|+r_2=\binom{r_1}{2}+r_2$ edges, let $G_2''$ be the colex graph with $|E(G_2)|-r_2=\binom{r_2-2}{2}+r_2-3$ edges. Then by Lemma~\ref{easy:convex}, 
	
	\begin{align*}
		k_t(G_1\cup G_2) &= \binom{r_1}{t}+\binom{r_2-2}{t}+\binom{r_2-1}{t-1}+\binom{r_2-2}{t-1}\\
		&< \binom{r_1}{t}+\binom{r_2}{t-1} +\binom{r_2-2}{t}+\binom{r_2-3}{t-1}\\
		&= k_3(G_1''\cup G_2''),
	\end{align*}
	a contradiction.
	
\end{proof}

\section{Proof of Theorem~\ref{cycle} and Corollary~\ref{edgecor}}
\label{sec:cor}

Both proofs follow from Theorem~\ref{edge}.

\begin{proof}[Proof of Theorem~\ref{cycle}]
Let $G$ be an $n$-vertex $C_{\ge k}$-free graph, where $n = q(k-2) + r$ with $1 \le r \le k-2$. By Theorem~\ref{kopylov}, the number of edges in $G$ is at most $q\binom{k-1}{2} + \binom{r}{2}$. Then, by Theorem~\ref{edge}, $k_t(G) \le k_t(qK_{k-1} \cup L)$, where $L \in \C_{t,k}(\binom{r}{2})$. Since $k_t(L)\le k_t(K_r)$, the upper bound in Theorem~\ref{cycle} follows immediately.
	
One can easily check that all the graphs described in Theorem~\ref{cycle} (i.e. graphs in $\mathcal{G}_{tree}(qK_{k-1}, K_r)$ and $\mathcal{G}'$) have indeed extremal structures. For the other direction, assume that $k_t(G)= q\binom{k-1}{t}+ \binom{r}{t}$. By Theorem~\ref{edge}, $G\in \Gf(qK_{k-1}, \C_{t, k}(\binom{r}{2})$. If $r \ge t$, then $K_r$ is the only graph in $\C_{t,k}(\binom{r}{2})$. Since $G$ has exactly $n = q(k-2) + r$ vertices it follows that $G \in \G(qK_{k-1}, K_r)$, as desired. If $r<t$, since $|V(G)|=n$, it follows immediately that $G\in \mathcal{G}'$. 
\end{proof}

\begin{proof}[Proof of Corollary~\ref{edgecor}]
Let $G$ be a $P_k$-free graph with $q\binom{k-1}{2}+\binom{r}{2}+s$ edges where $0\le s<r< k-1$. Then $G$ is also a $C_{\ge k}$-free graph. Thus, Theorem~\ref{edge} provides the desired upper bound on $k_t(G)$. If the upper bound is tight, i.e. $k_t(G)=q\binom{k-1}{t}+\binom{r}{k-1} +\binom{s}{t-2}$, then $G\in \mathcal{G}_{forest} (qK_{k-1}, L)$ for some $L\in \mathcal{C}_{t, k}(\binom{r}{2}+s)$. Since $G$ is $P_k$-free, any cliques of order $k-1$ are disconnected from the rest of the graph. Then $G$ is isomorphic to $qK_{k-1}\cup L$ and $L\in \mathcal{L}_{t, k}(\binom{r}{2}+s)$. Thus, $G$ is one of the extremal graphs. 
\end{proof} 

\section{Concluding remarks}

As discussed in \cite{AS,GMP}, an interesting question is to consider the generalized Erd\H{o}s-S\'os conjecture where one wishes to maximize the number of $K_t$'s in an $n$-vertex graph that does not contain a copy of a fixed tree $T$. Asymptotically tight bounds were proven in \cite{GMP} on this number assuming that the original Erd\H{o}s-S\'os conjecture \cite{E} holds for every tree. Theorem~\ref{path} and a result by Chase \cite{C}, for example, are the special cases of this generalized problem where we avoid the path $P_k$ and the star $K_{1, \Delta}$ respectively. For a related result on $\ex(n,H,T)$ where $H$ is a general graph and $T$ is a tree, curious readers can refer to \cite{L}.

Similar to Corollary~\ref{edgecor}, the variation where the number of edges is fixed instead of the number of vertices when excluding other trees can also be considered. Thus, a natural objective is to find some non-trivial upper bound for this edge variant, even for certain small classes of forbidden trees. 

In an alternative direction, one might consider the following generalization of Theorem~\ref{path} and Corollary~\ref{edgecor}: for $0 < s < t$, $m$, and $k$ positive integers, determine the maximum number of $K_t$'s in a $P_k$-free graph with $m$ copies of $K_s$.  Once again, one may replace the forbidden structure with any other tree and consider similar interesting questions. 

\section{Acknowledgement}
We thank D\'aniel Gerbner for bringing the reference \cite{GMP} to our attention and extend our appreciation to the anonymous reviewers for their careful reading.

\end{document}